# Nonlinear MPC with PWM applied on a small satellite

Jinaykumar Patel[1]
*The University of Texas at Arlington, Arlington, TX, 76010, USA*

This paper employs Pulse Width Modulation (PWM) method to a Nonlinear Model Predictive Control (NMPC) for satellite attitude control and detumbling. Magnetic torquers have been used as actuators for the stabilization and attitude control of small satellites. NMPC generates the control input which is continuous and smooth. As a result, it is challenging for the on-board actuators to produce these continuous magnetic moments as a control input. Thus, PWM is applied to discretize the control input reducing the load on the actuators. Simulations are presented for detumbling and attitude control to illustrate the effectiveness and feasibility of the proposed method.

## Nomenclature

| | | |
|---|---|---|
| $i$ | = | inclination |
| $J$ | = | cost function |
| $I_B$ | = | moment of inertia in body frame |
| $m$ | = | magnetic dipole moment |
| $u_{max}$ | = | maximum control input |
| $q$ | = | quaternion |
| $\omega_B$ | = | Angular velocity in body frame |
| $Q$ and $R$ | = | weight matrices |
| $\tau$ | = | torque generated by magnetic torquers |
| $T_s$ | = | sampling time |
| $k$ | = | $k^{th}$ time step of Model Predictive Control |
| $x_0$ | = | Initial state vector |
| $x_{ref}$ | = | Reference or Final state vector |
| $B$ | = | Earth's magnetic field |
| p | = | prediction horizon |

## I. Introduction

The stabilization of a satellite after deployment is always an important aspect of every mission. Satellite may get lost during this phase if it is not stabilized. The angular rates at the deployment are large than those desired for the attitude maneuvering. Therefore, a detumbling control algorithm must be implemented to achieve stability in minimum time. Once the angular rate drops, the controller can be used for attitude control and maneuvering. Magnetic actuators have been widely used for stabilization and attitude control of small satellites due to restrictions on size and weight [1]. Magnetic torquers have commonly been applied to spacecraft to attenuate angular rates [2]. They have also been extensively used with momentum wheels and for damping of spin satellites [2]. Basically, magnetic torquers creates a magnetic dipole moment which interacts with the Earth's magnetic field and thus, a torque is generated. Magnetic dipole moment act as a control input to the system. However, magnetic actuators require dense magnetic flux, which makes the control system under-actuated [3], [4]. Nevertheless, if satellite orbit inclination with equator is large, the magnetic field vector periodically changes and as a result, it is possible to apply torques in any direction [5].

---
[1] Graduate Student, Department of Mechanical and Aerospace Engineering.

As the magnetically actuated system is underactuated, the design of controllers for such a system is therefore a pretty challenging task. In this work, feedback control method, nonlinear MPC is employed for the detumbling and attitude control of satellite. MPC optimizes a model to generate series of control inputs that minimizes the cost function over a finite control horizon while subjected to numerous constraints. Future plant outputs and control actions are predicted based on prior values. The advantages of using nonlinear MPC is that it can incorporate nonlinear constraints and also, can incorporate multiple constraints simultaneously. Moreover, MPC is a closed loop control and its cost function includes the control input. The MPC finds control inputs which are continuous and smooth, which creates burden on the control actuators. Therefore, PWM technique is employed [6], which discretizes the smooth control inputs to a discrete value. This discretization will help reducing the burden on the actuators.

This paper demonstrates the feasibility of the nonlinear MPC applied with the PWM for the detumble and attitude control of small satellite. Sections II includes the rotational kinematics and dynamics equation of motion. Section III and IV describes the MPC and PWM problem formulation, respectively. The simulation section V shows the detumbling and attitude control simulations. In detumbling, the satellite's angular rates are attenuated, while in attitude control simulation, a desired attitude is attained. Finally, section VI and VII concludes the paper.

## II. Equations of Motion

All equations are derived assuming a body-fixed frame with origin at the center of mass of the satellite and the principal axes of inertial aligned with the axis of frame.

### A. Rotational Kinematics

The rotational kinematics equations in quaternion form are given as [7]:

$$\begin{bmatrix} \dot{q}_1 \\ \dot{q}_2 \\ \dot{q}_3 \end{bmatrix} = \frac{1}{2} S(\omega_B) \begin{bmatrix} q_1 \\ q_2 \\ q_3 \end{bmatrix} + \frac{1}{2} q_4 \omega_B \tag{1}$$

$$\dot{q}_4 = -\frac{1}{2} \omega_B^T \begin{bmatrix} q_1 \\ q_2 \\ q_3 \end{bmatrix} \tag{2}$$

Combining equations Eq. (1) and (2), we can write as:

$$\dot{\boldsymbol{q}} = M(\boldsymbol{q}) \, \boldsymbol{\omega}_B \tag{3}$$

where, $\boldsymbol{q} = \begin{bmatrix} q_1 \\ q_2 \\ q_3 \\ q_4 \end{bmatrix}$, $\boldsymbol{\omega}_B = \begin{bmatrix} \omega_x \\ \omega_y \\ \omega_z \end{bmatrix}$, $S(\omega_B) = \begin{bmatrix} 0 & \omega_z & -\omega_y \\ -\omega_z & 0 & \omega_x \\ \omega_y & -\omega_x & 0 \end{bmatrix}$ and $M(\boldsymbol{q}) = \frac{1}{2} \begin{bmatrix} q_4 & -q_3 & q_2 \\ q_3 & q_4 & -q_1 \\ -q_2 & q_1 & q_4 \\ -q_1 & -q_2 & -q_3 \end{bmatrix}$

### B. Rotational Dynamics

Rotational dynamics is modelled assuming the satellite as a rigid body on which the only torque acts is due to magnetic dipole. This is the torque due to the interaction of the magnetic dipole moment produced by magnetic torques with the Earth's magnetic field. The rotational dynamics motion in the body-fixed frame is given by the Euler's equation [7]:

$$\boldsymbol{I}_B \dot{\boldsymbol{\omega}}_B = S(\boldsymbol{\omega}_B) \boldsymbol{I}_B \boldsymbol{\omega}_B + \boldsymbol{\tau}_B \tag{4}$$

The torque generated through the interaction of magnetic torquers with the Earth's magnetic field is given as:

$$\boldsymbol{\tau}_B = \boldsymbol{m} \times \boldsymbol{B}(t) \tag{5}$$

where, $\boldsymbol{m} = [m_x, m_y, m_z]^T$ is the magnetic dipole moment and $\boldsymbol{B}(t) = [B_x, B_y, B_z]^T$ is the magnetic field of the earth.

Eq. (5) can further be extended as

$$\boldsymbol{\tau}_B = \boldsymbol{m} \times \boldsymbol{B}(t) = \begin{bmatrix} m_y B_z - m_z B_y \\ m_z B_x - m_x B_z \\ m_x B_y - m_y B_x \end{bmatrix} \qquad (6)$$

Combining Eq. (4) and Eq. (6), we can arrive at the $\dot{\boldsymbol{\omega}}_B$ expression:

$$\begin{bmatrix} \dot{\omega}_x \\ \dot{\omega}_y \\ \dot{\omega}_z \end{bmatrix} = \begin{bmatrix} \frac{1}{I_x}\{(I_y - I_z)\omega_y\omega_z + m_y B_z - m_z B_y\} \\ \frac{1}{I_y}\{(I_z - I_x)\omega_z\omega_x + m_z B_x - m_x B_z\} \\ \frac{1}{I_z}\{(I_x - I_y)\omega_x\omega_y + m_x B_y - m_y B_x\} \end{bmatrix} \qquad (7)$$

### III. Model Predictive Control Formulation

The basic idea of the Model Predictive Control is to numerically optimize a model to obtain a sequence of control input that minimizes a cost (or objective) function over a finite receding horizon, while subjected to certain constraints. The control actions are periodically recomputed at each sampling instant with the current state estimate as an initial condition, thereby providing a feedback action. The main goal of the MPC is to reduce the errors between the predicted states ($x$) and the reference state ($x_{ref}$). Thus, the cost function, denoted as $J$, defined as follows:

$$J = \int_t^{t+T} \left\{ \left(\boldsymbol{x}(k) - \boldsymbol{x}_{ref}\right)^T \boldsymbol{Q} \left(\boldsymbol{x}(k) - \boldsymbol{x}_{ref}\right) + \boldsymbol{u}(k)^T \boldsymbol{R}\, \boldsymbol{u}(k) \right\} dk \qquad (8)$$

where, $\boldsymbol{x} = [q_1, q_2, q_3, q_4, \omega_x, \omega_y, \omega_z]^T$, $\boldsymbol{u} = [u_x, u_y, u_z]^T$, $\boldsymbol{Q}$ and $\boldsymbol{R}$ are symmetric positive definite matrices which penalizes the state values and the control input, respectively.

Thus, the optimal control problem is to minimize cost function, $J$ subject to the following constraints:

$$\min_{\boldsymbol{u}} J$$

$$\text{subject to} \begin{cases} \dot{\boldsymbol{q}} = \frac{1}{2} \begin{bmatrix} q_4 & -q_3 & q_2 \\ q_3 & q_4 & -q_1 \\ -q_2 & q_1 & q_4 \\ -q_1 & -q_2 & -q_3 \end{bmatrix} \boldsymbol{\omega}_B \\ \boldsymbol{I}_B \dot{\boldsymbol{\omega}}_B = S(\boldsymbol{\omega}_B)\boldsymbol{I}_B \boldsymbol{\omega}_B + \boldsymbol{\tau}_B \\ \boldsymbol{\tau}_B = \boldsymbol{m} \times \boldsymbol{B}(t) = \begin{bmatrix} m_y B_z - m_z B_y \\ m_z B_x - m_x B_z \\ m_x B_y - m_y B_x \end{bmatrix} \\ |m_i| \leq u_{max} \quad \forall\ i = x, y, z \end{cases}$$

### IV. Pulse Width Modulation

A pulse width modulation technique is utilized to convert continuous and smooth control input signal into a discrete signal which will help in reducing the load on the actuator. An algorithm presented in [8] and [9], is employed here, where previous value of input is used to determine the input at the current time stamp and the algorithm tends to keep the control input value same as the previous value. The algorithm uses the logic depending on the previous input value and sort the signal into discrete fraction values of the maximum control input ($u_{max}$). An example of branching of algorithm is as follows:

$$u = \begin{cases} u_{max}, & \text{if } u_{max} \geq \frac{2}{3} \\ \frac{2}{3}u_{max}, & \text{if } \frac{2}{3}u_{max} > u_{mpc} \geq \frac{1}{3} \\ \frac{1}{3}u_{max}, & \text{if } \frac{1}{3}u_{max} > u_{mpc} \geq 0 \\ 0, & \text{if } 0 > u_{mpc} \geq -\frac{1}{3}u_{max} \\ -\frac{1}{3}u_{max}, & \text{if } -\frac{1}{3}u_{max} > u_{mpc} \geq -\frac{2}{3}u_{max} \\ -\frac{2}{3}u_{max}, & \text{if } -\frac{2}{3}u_{max} > u_{mpc} \geq -u_{max} \\ -u_{max}, & \text{if } -u_{max} > u_{mpc} \end{cases}$$

## V. Simulation and Results

In this part, simulations of satellite detumbling and attitude control are performed where the feasibility of the proposed model predictive control algorithm with PWM is showcased. A sun-synchronous orbit [9] is used for the subsequent simulations in this section. The orbital parameters for the sun-synchronous orbit are described in Table 1.

**Table 1 Orbital Elements**

| Parameters | Values |
| --- | --- |
| Semi-major axis | 6691.6 km |
| Inclination | 96.7° |
| Eccentricity | 0.046440 |
| Mean Anomaly | 240.49° |
| Right Ascension of Ascending Node | 100.90° |
| Argument of Perigee | 119.70° |

### A. Magnetic Field Model

As from Eq. 6, we need the Earth's magnetic field for the calculation of torque. There are various models developed for the Earth's magnetic field like, International Geomagnetic Reference Field (IGRF) model [10]. However, IGRF model involves harmonic coefficients and to use this model for on-board computation is expensive. Hence, we employ an approximated model as described in [5]. The components of magnetic field $[B_x, B_y, B_z]$ in the orbital frame of reference are describe by following equations:

$$\begin{bmatrix} B_x \\ B_y \\ B_z \end{bmatrix} = D_m \begin{bmatrix} \frac{3}{2}\sin i \sin 2\eta \\ -\frac{3}{2}\sin i (\cos 2\eta - \frac{1}{3}) \\ -\cos i \end{bmatrix} \quad (9)$$

where $\eta = \theta + \omega_e$, $D_m = -\frac{M_e}{r^3}$, $M_e = 8.1 \times 10^{25}\ gauss.cm^3$, $i$ is the inclination, $\omega_e$ is the argument of perigee, $\theta$ is true anomaly, and $r$ is the distance of the satellite from the center of the earth.
Using the orbital elements from the Table 1, the time varying magnetic field is shown in Figure 1.

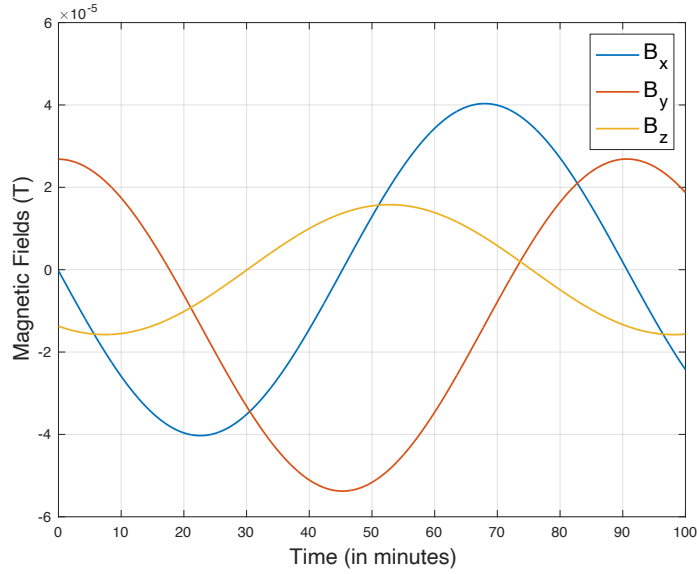

**Figure 1** Magnetic Field for sun-synchronous orbit

**Table 2** MPC parameters for detumble

| Parameters | Values |
| --- | --- |
| Control Input, $u$ | $[-0.1, 0.1]\ A.m^2$ |
| Sampling time, $T_s$ | $2\ sec$ |
| Q | diag ([0, 0, 0, 0, 500 1000 250]) |
| R | diag ($[10^{-8}\ \ 10^{-8}\ \ 10^{-8}]$) |
| Initial State $x_0$ | $[0, 0, 0, 1, 4(deg/s), 3(deg/s), -3(deg/s)]$ |
| Reference or Final State $x_f$ | $[0, 0, 0, 1, 0, 0, 0]$ |
| Prediction Horizon | 10 |

**Table 3** Moment of Inertia of satellite

| **Moment of inertia** | **Value** |
| --- | --- |
| $I_{xx}$ | $0.020\ kg.m^2$ |
| $I_{yy}$ | $0.030\ kg.m^2$ |
| $I_{zz}$ | $0.040\ kg.m^2$ |

## B. Detumbling Simulation

In this section, detumblimg of satellite's angular velocities is demonstrated. Table 2 describes the parameters for the simulation and the moment of inertia of satellite used for the simulation are described in Table 3. Figure 2 illustrates the detumbling result. Starting from initial state $x_0$, the model achieve stability within 25 minutes, i.e. it reaches the final state $x_f$. Initially there is large variations in angular velocity but ultimately it converges close to 0.

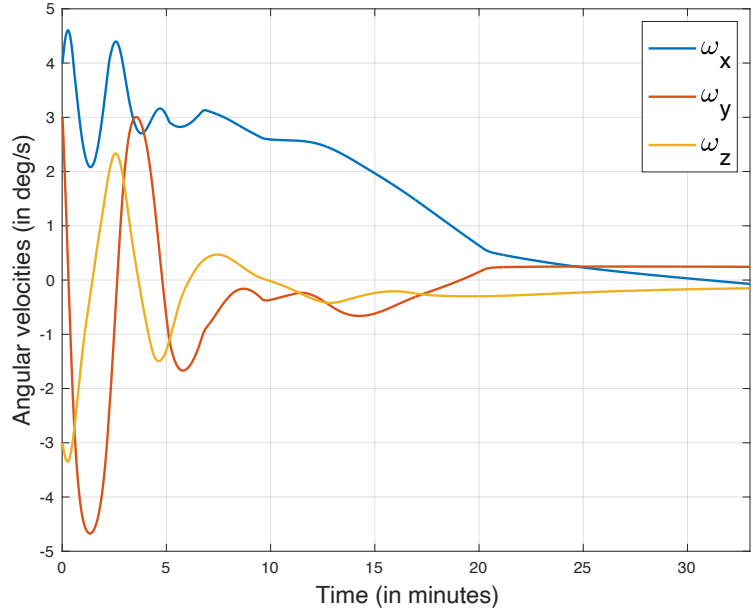

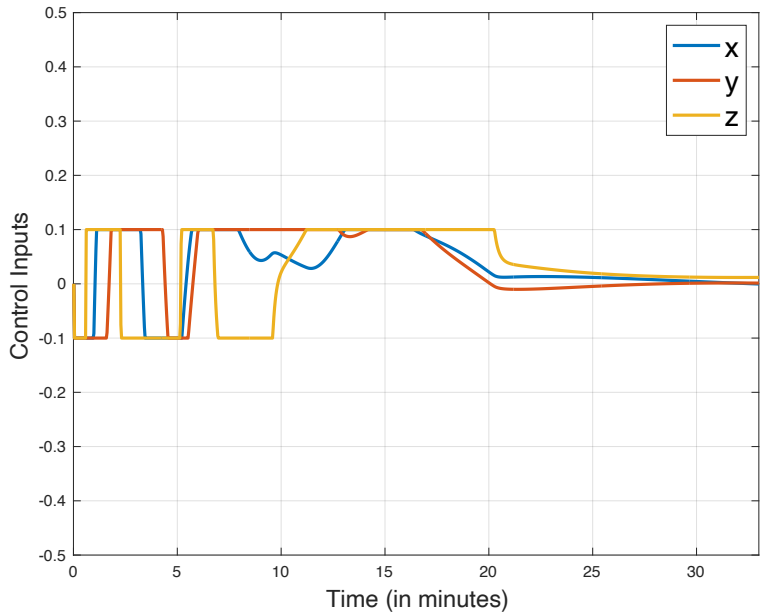

**Figure 2   (a) Angular velocity and (b) control input time history**

### C.  Attitude Control Simulation

This section demonstrates attitude control simulation results. Table 4 describes the moment of inertia used for the simulation and Table 4 highlights the MPC model parameters used for attitude control simulation, assuming the initial angular velocity components to be 0.  Figure 3 illustrates that the attitude control result. Clearly, from Figure 3 (a), starting from initial state $x_0$, the model achieve stability within 25 minutes, i.e. it reaches the final state $x_f$.

Table 4  MPC parameters for attitude control

| Parameters | Values |
| --- | --- |
| Control Input, $u$ | $[-0.1, 0.1]\ A.m^2$ |
| Sampling time, $T_s$ | $30\ sec$ |
| Q | diag $([20, 20, 20, 20, 10^4\ 10^4\ 10^4])$ |
| R | diag $([10^{-8}\ \ 10^{-8}\ \ 10^{-8}\ ])$ |
| Initial State $x_0$ | $[0, 0.1, 0, 1, 0, 0, 0]$ |
| Reference or Final State $x_f$ | $[0, 1, 0, 0, 0, 0, 0]$ |
| Prediction Horizon | 10 |

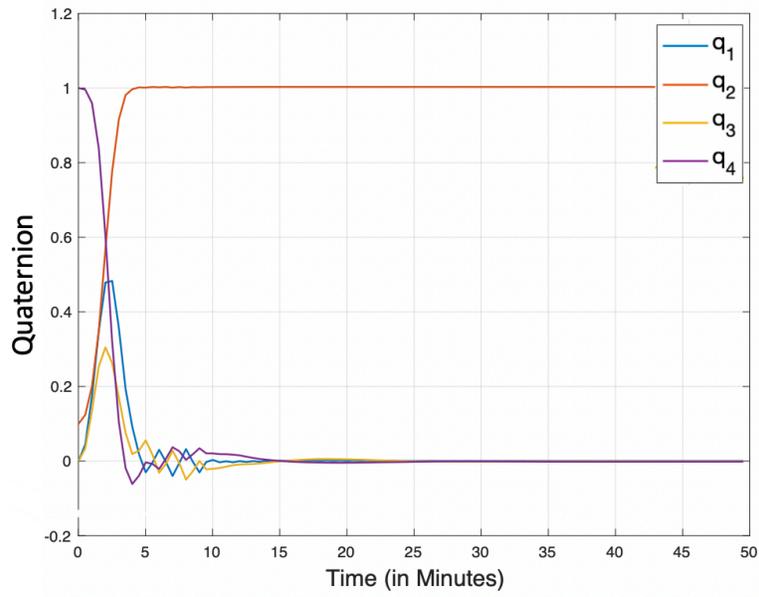

(a)

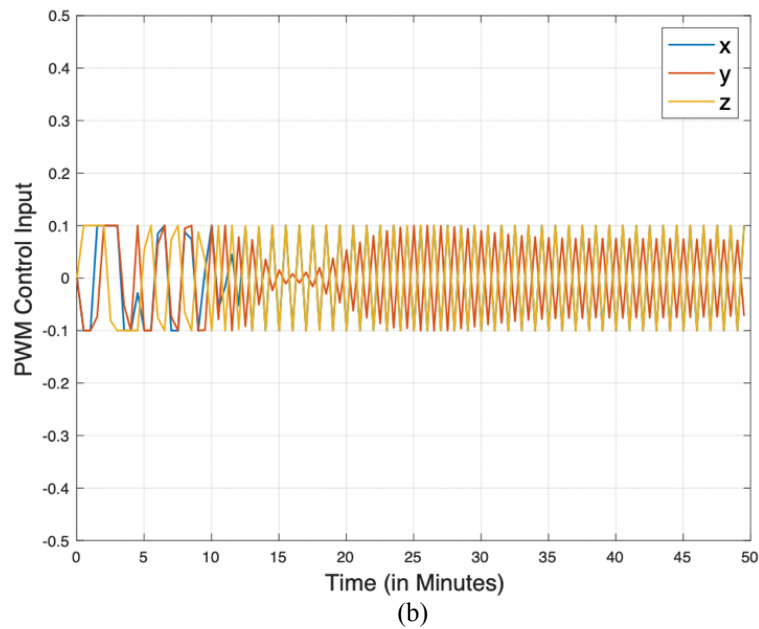

(b)

**Figure 3  (a) Quaternion and (b) control input time history**

## D. Comparison of control input (With PWM and Without PWM)

As shown in Figure 4, PWM discretizes the continuous control input. Figure 4(a) shows continuous control inputs from MPC controllers and Figure 4(b) presents the discrete input from PWM method. We can clearly see that the with PWM, the control input is discretized into step values.

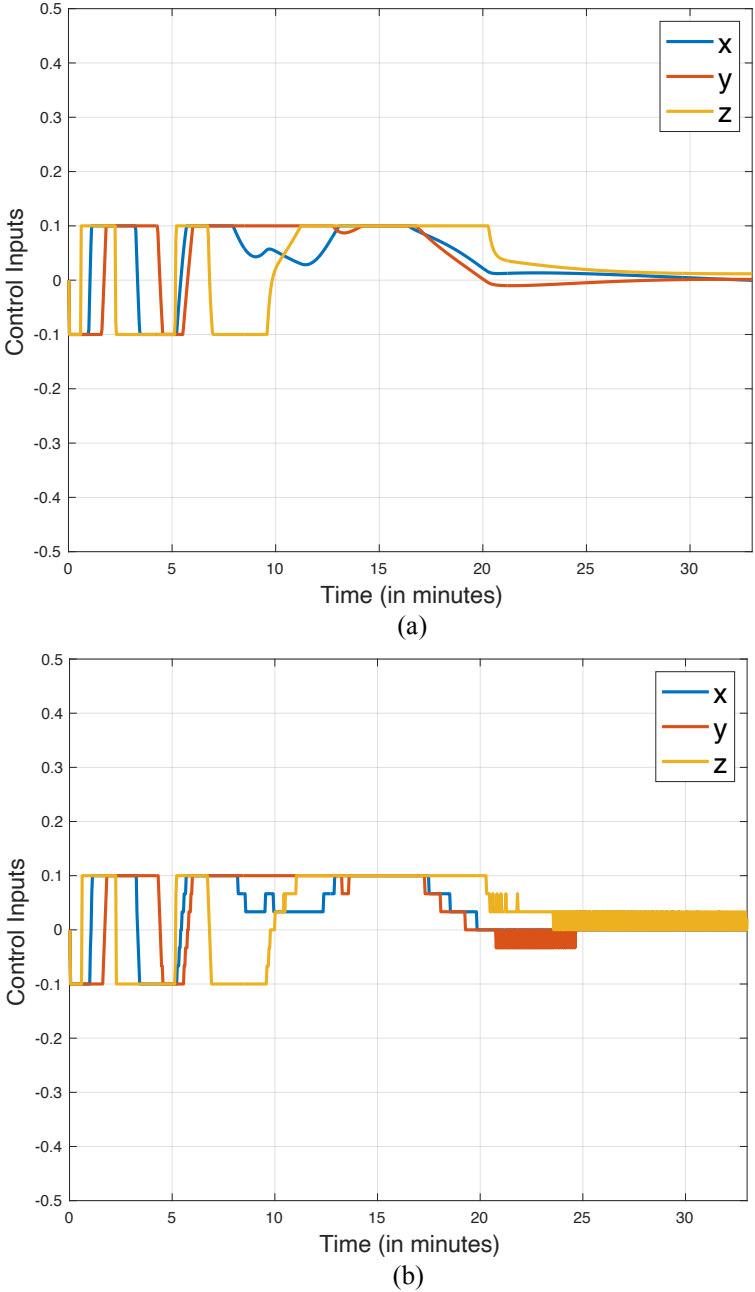

**Figure 4    (a) Control Input without using PWM and (b) with PWM**

## VI. Conclusion

In conclusion, this paper has implemented the PWM method to convert the smooth MPC control input signals into discrete values. The discretization reduces the burden on the control actuators. The design and numerical simulation of detumbling and angular rate was presented using magnetorquers. Nonlinear MPC with PWM is able to detumble the small satellite and also, able to control the attitude of satellite using magnetic torquers. The proposed MPC strategy has shown that detumbling is achieved within allocated time. Thus, the introduction of PWM scheme in conjunction with a closed loop MPC satisfies all control objectives in the simulations.

## VII. Future Work

Future work includes using this model for different orbital parameters, i.e. for different orbits of satellite. Also, further analysis on the stability will help improve the current model. Moreover, this analysis can even be extended to spacecraft.

## References

[1] D. K. Giri and M. Sinha, "Fast Terminal Sliding-Mode Fault-Tolerant Attitude Control of Magnetically Actuated Satellite," *J. Spacecr. Rockets*, vol. 56, no. 5, pp. 1636–1645, Sep. 2019, doi: 10.2514/1.A34475.
[2] K. Kondo, I. Kolmanovsky, Y. Yoshimura, M. Bando, S. Nagasaki, and T. Hanada, "Nonlinear Model Predictive Detumbling of Small Satellites with a Single-Axis Magnetorquer," *J. Guid. Control Dyn.*, vol. 44, no. 6, pp. 1211–1218, Jun. 2021, doi: 10.2514/1.G005877.
[3] X. Huang and Y. Yan, "Fully Actuated Spacecraft Attitude Control via the Hybrid Magnetocoulombic and Magnetic Torques," *J. Guid. Control Dyn.*, vol. 40, no. 12, pp. 3358–3360, Dec. 2017, doi: 10.2514/1.G002925.
[4] M. Wood, W.-H. Chen, and D. Fertin, "Model predictive control of low earth orbiting spacecraft with magneto-torquers," in *2006 IEEE Conference on Computer Aided Control System Design, 2006 IEEE International Conference on Control Applications, 2006 IEEE International Symposium on Intelligent Control*, Munich, Germany, Oct. 2006, pp. 2908–2913. doi: 10.1109/CACSD-CCA-ISIC.2006.4777100.
[5] J. Cubas, A. Farrahi, and S. Pindado, "Magnetic Attitude Control for Satellites in Polar or Sun-Synchronous Orbits," *J. Guid. Control Dyn.*, vol. 38, no. 10, pp. 1947–1958, Oct. 2015, doi: 10.2514/1.G000751.
[6] J. Holtz, "Pulsewidth modulation-a survey," *IEEE Trans. Ind. Electron.*, vol. 39, no. 5, pp. 410–420, Oct. 1992, doi: 10.1109/41.161472.
[7] D. C. Elliott, "AE 5302 - ADVANCED FLIGHT MECHANICS LECTURE NOTES," p. 632.
[8] H. Seguchi and T. Ohtsuka, "Nonlinear receding horizon control of an underactuated hovercraft," *Int. J. Robust Nonlinear Control*, vol. 13, no. 3–4, pp. 381–398, Mar. 2003, doi: 10.1002/rnc.824.
[9] K. Kondo, Y. Yoshimura, S. Nagasaki, and T. Hanada, "Pulse Width Modulation Method Applied to Nonlinear Model Predictive Control on an Under-actuated Small Satellite," presented at the AIAA Scitech 2021 Forum, VIRTUAL EVENT, Jan. 2021. doi: 10.2514/6.2021-0559.
[10] P. Alken *et al.*, "International Geomagnetic Reference Field: the thirteenth generation," *Earth Planets Space*, vol. 73, no. 1, p. 49, Dec. 2021, doi: 10.1186/s40623-020-01288-x.